\documentclass[12pt,a4paper]{article}
\usepackage[english]{babel}
\usepackage{indentfirst}
\usepackage[dvips]{graphicx}
\usepackage{amstext}
\usepackage{amsfonts}
\usepackage{textcomp}
\usepackage{amssymb}
\usepackage{amscd}
\usepackage[dvips]{graphicx}
\usepackage{setspace}
\newcommand {\demo}{{\it \bf Proof:  }}
\newcommand {\fim}{$\rightline{$\square$}$\vskip 2pc}
\newcommand {\nl}{\newline}
\newcommand {\cl}{\centerline}

\newcommand {\N}{\mathbb{N}}

\newcommand {\F}{\mathbb{F}}
\newcommand {\G}{\mathbb{G}}

\newcommand {\Lh}{\mathcal{L}}
\newcommand {\La}{\Lambda}
\newcommand {\al}{\alpha}
\newcommand {\til}{\widetilde}
\newcommand {\lb}{\linebreak}
\newcommand {\esp}{\hskip 1pc}

\newtheorem{teorema}{Theorem}
\newtheorem{lema}{Lemma}
\newtheorem{definicao}{Definition}

\title{Markov subshifts and partial representation of $\F_n$}
\author{Danilo Royer}

\begin{document}
\maketitle

\begin{abstract} In this paper we fix a set $\Lambda^*$ of positive
elements of the free group $\F_n$ (e.g.~the set of finite words
occurring in a Markov subshift) as well as $n$ partial isometries on a
Hilbert space $H$.  Based on these we define a map $S:\F_n\rightarrow
\Lh(H)$ which we prove to be a partial representation of $\F_n$ on $H$
under certain conditions studied by Matsumoto.

\end{abstract}

\section*{Introduction}

Considering a Markov subshift on an alphabet $\{g_1,\ldots, g_n\}$,
R. Exel proved in $\cite{Exel}$ that $n$ partial isometries on a Hilbert 
space
$H$, satisfying the corresponding Cuntz--Krieger relations, give rise
to a partial representation of the free group $\F_n$ on $H$, that is,
a map $S:\F_n\longrightarrow\Lh(H)$, satisfying $S(t^{-1})=S(t)^*$ and
$S(tr)S(r^{-1})=S(t)S(r)S(r^{-1})$ for all $r,t$ in $\F_n$.

In this work we fix a set $\La^*$ of positive elements of $\F_n$
which, among other requirements is assumed to be closed under
sub-words, and we take a set $\{S_1,\ldots,S_n\}$ of partial
isometries on $H$. We define a map $S:\F_n\longrightarrow\Lh(H)$ by
$S(r_1\ldots r_k)=S(r_1)\ldots S(r_k)$, where $S(r_i)=S_j$ if
$r_i=g_j$, $S(r_i)=S_j^*$ if $r_i=g_j^{-1}$ and $r=r_1\ldots r_k$ is in  
reduced form.

Under certain conditions studied by Matsumoto in $\cite{Matsumoto}$, we 
prove that
the map $S$ is a partial representation of $\F_n$ on $H$.  Since
Matsumoto's conditions generalize the Cuntz--Krieger relations our
result is a generalization of Exel's result mentioned above.

This paper is based on the author's Master's thesis at the Federal
University of Santa Catarina under the supervision of Ruy Exel.

\section*{Partial Representations of $\F_n$}
\par Let us consider the Free Group $\F_n$ generated by a set of $n$ elements, $G=\{g_1,\ldots
,g_n\}$. The elements of $\F_n$ can be written in the form $r=r_1\ldots  r_k$
where each $r_i\in G\cup G^{-1}$. We say that $r$ is in reduced form if
$r_i\neq r_{i+1}^{-1}$, for each $i$. Two elements $r=r_1\ldots r_k$ and
$s=s_1\ldots s_l$ of $\F_n$, in reduced form, are equal if and only if $l=k$
and $r_i=s_i$, for all $i$. In this way, each element, in reduced form, have
unique representation and we define its length by the number of components,
that is, if $r=r_1\ldots r_k$ is in reduced form then $r$ have length $k$, which
will be denoted by $|r|=k$. A element $r=r_1\ldots r_k$ of $\F_n$, in reduced
form, is called a positive element if $r_i\in G,$ for all $i$, and the set of
all positive elements will be called $P$. We consider $e$ a element of $P$.
\par Let us fix a set $\La^*\subseteq P$ with the following properties:
\begin{enumerate}
\item [$\bullet$] $e\in \La^*$,
\item [$\bullet$] $G=\{g_1,\ldots,g_n\}\subseteq \La^*$,
\item [$\bullet$] $\La^*$ is closed under sub-words, that is, if $\nu=\nu_1\ldots \nu_k\in \La^*$ then each element
of the form $\nu_i\ldots\nu_{i+j}$ with $i=1\ldots k, j\in \N$ is a element of
$\La^*$.

\end{enumerate}

\par For all $\mu\in \La^*$ we define the following sets:
$$L_{\mu}^1=\{g_j\in G| j=1,\ldots ,n,\esp \mu g_j\notin \La^*\},$$
$$L_{\mu}^k=\{\nu=\nu_1\ldots \nu_k\in \La^* | \mu \nu_1\ldots \nu_{k-1}\in
\La^*, \mu \nu\notin \La^*\},\hskip 2pc \forall k\in \N.$$

\begin{lema}\label{lt1}
Let $\mu\in\La^*$ and $r,s\in P$. If $vr=v's$, where $v\in L_\mu^k$ and $v'\in
L_\mu^l$, then $v=v'$.
\end{lema}
\demo Suppose by contradiction that $v\neq v'$. Then $|v|\neq|v'|$, because
otherwise, $v_1\ldots v_kr=v_1'\ldots v_k's$, from where it follows that
$v=v'$. Without loss of generality suppose $|v|>l$, write $v=v_1\ldots
v_l\ldots v_k$ and $v'=v_1'\ldots v_l'$. Since \lb $v_1\ldots v_l\ldots
v_kr=vr=v's=v_1'\ldots v_l's$, then $v_1\ldots v_l=v_1'\ldots v_l'$, and
therefore \lb  $v=v'v_{l+1}\ldots v_k$. Since $v'\in L_\mu^l$, by definition of
$L_\mu^l$, $\mu v' \notin \La^*$, hence \lb $\mu v_1\ldots v_{k-1}=\mu v'
v_{l+1}\ldots v_{k-1}\notin \La^*$. That is a contradiction, because $v\in
L_\mu^k$ and so $v=v'$.\nl \fim 

\par Let us consider a Hilbert space $H$ and a set of partial
isometries $\{S_1,\ldots ,S_n\}\subseteq\Lh(H)$. Recall that $S_i$ is a partial
isometry if $S_iS_i^*S_i=S_i$. Define a map \nl \cl{$S:\F_n\longrightarrow
\mathcal{L}(H)$}\nl \cl{$r=r_1\ldots r_k\mapsto S(r_1)\ldots S(r_k)$}\nl where
$r$ is in reduced form, $S(r_i)=S_j$ if $r_i=g_j$ and $S(r_i)=S_j^*$ if
$r_i=g_j^{-1}$. By convention, $S(e)=I$, where $I$ is the identity operator on
$H$. In this way, for all $r\in \F_n$ we have an operator $S(r)\in
\mathcal{L}(H)$.This operator will also be called $S_r$. We will suppose that
our set of partial isometries $\{S_1,\ldots,S_n\}\subseteq \Lh (H)$ generated a
map $S$ which satisfies:
\begin{enumerate}
\item [$(M_1)$] $\sum\limits_{i=1}^n S_iS_i^*=I$;
\item [$(M_2)$] For all $\mu$ and $\nu$ in $\La^*$ the operators
$S_\mu S_\mu^*$ and $S_\nu^* S_\nu$ commute;
\item [$(M_3)$] $I-S_i^*S_i=\sum\limits_{k=1}^\infty\sum\limits_{\nu\in
L_i^k}S_\nu S_\nu^*$,  $i=1,\ldots ,n$.
\end{enumerate}

Note that for all $i$, $S_iS_i^*$ is idempotent and self-adjoint, and so a
projection. By $(M_1)$, $\sum\limits_{i=1}^n S_iS_i^*$ is a projection and
therefore $S_iS_i^*$ and $S_jS_j^*$ are orthogonal, for all $i\neq j$. So
$$S_i^*S_j=(S_i^*S_iS_i^*)(S_jS_j^*S_j)=S_i^*(S_iS_i^*S_jS_j^*)S_j=0$$ whenever
$i\neq j$.

\begin{lema}\label{lt2}
For all $\mu\in\La^*$, $S_\mu=S_\mu S_\mu^* S_\mu$.
\end{lema}
\demo The proof will be by induction on $|\mu|$. For $|\mu|=1$, $S_\mu=S_\mu
S_\mu^* S_\mu$ by hypothesis. Suppose $S_\mu=S_\mu S_\mu^* S_\mu$ for all
$\mu\in \La^*$ with $|\mu|=k$, and consider $\nu\in\La^*$, with $|\nu|=k+1$.
Then $\nu=\alpha g_j$, with $|\alpha|=k$, and $$S_\nu S_\nu^* S_\nu=S_{\alpha
g_j}S_{\alpha g_j}^* S_{\alpha g_j}=S_\alpha S_{g_j}S_{g_j}^*S_\alpha^*S_\alpha
S_{g_j}=S_\alpha S_\alpha^*S_\alpha S_{g_j}S_{g_j}^*S_{g_j}=S_\alpha
S_{g_j}=S_\nu.$$ \fim

\begin{lema}\label{lt3}
Let $\alpha\in P$ and $\nu\in\La^*$.
\begin{enumerate}
\item [a)] If $|\al|\geq|\nu|$ then $S_\nu S_\nu^* S_\al=\left\{\begin{array}{cc}
  S_\al & \text{ if }\al=\nu r \text{ for some }r\in P\\
  0 & \text{ otherwise}
\end{array}\right .$

\item [b)] If $|\al|<|\nu|$ then $S_\nu S_\nu^* S_\al=\left\{\begin{array}{cc}
  S_\nu S_r^*& \text{ if }\nu=\al r \text{ for some }r\in P\\
  0 & \text{ otherwise}
\end{array}\right .$
\end{enumerate}
\end{lema}
\demo
\begin{enumerate}
\item [a)] Supposing that there exists $r$ in $P$ such that $\al=\nu r$, we have $$S_\nu S_\nu^* S_\al =
S_\nu S_\nu^* S_{\nu r}=S_\nu S_\nu^* S_\nu S_r=S_\nu S_r=S_\al.$$ On the other
hand, if $\alpha\neq\nu r$ for all $r\in P$, write $\al=\al_1\ldots \al_l\ldots
\al_k$, $\nu=\nu_1\ldots \nu_l$ and take the smallest index $i$ such that
$\al_i\neq\nu_i$. Then we have $\al_1\ldots \al_{i-1}=\nu_1\ldots \nu_{i-1}$,
and so $$S_\nu S_\nu^*S_\al=S_{\nu_1\ldots \nu_{i-1}\nu_i\ldots
\nu_l}S_{\nu_1\ldots \nu_{i-1}\nu_i\ldots \nu_l}^* S_{\al_1\ldots
\al_{i-1}\al_i\ldots \al_k}=$$ $$=S_{\nu_1\ldots \nu_{i-1}}S_{\nu_i\ldots
\nu_l}S_{\nu_i\ldots \nu_l}^*S_{\nu_1\ldots \nu_{i-1}}^* S_{\nu_1\ldots
\nu_{i-1}}S_{\al_i\ldots \al_k}=$$ $$=S_{\nu_1\ldots \nu_{i-1}}S_{\nu_1\ldots
\nu_{i-1}}^* S_{\nu_1\ldots \nu_{i-1}}S_{\nu_i\ldots \nu_l}S_{\nu_i\ldots
\nu_l}^*S_{\al_i\ldots \al_k}=0$$ because $S_{\nu_i}^*S_{\al_i}=0$.

\item [b)] Suppose $\nu=\al r$ for some $r\in P$. Then $$S_\nu S_\nu^* S_\al = S_{\al r}
S_{\al r}^* S_\al= S_\al S_r S_r^*S_\al^* S_\al= S_\al S_\al^* S_\al S_r S_r^*=
S_\al S_r S_r^*=S_{\al r}S_r^*=S_\nu S_r^*.$$ If $\nu\neq \alpha r$, for all
$r\in P$ as in (a), take the smallest index $i$ such that $\nu_i\neq\al_i$.
Then $\nu_1\ldots \nu_{i-1}=\al_1\ldots \al_{i-1}$ and $$S_\nu S_\nu^* S_\al=
S_{\nu_1\ldots \nu_{i-1}\nu_i\ldots \nu_k}S_{\nu_1\ldots \nu_{i-1}\nu_i\ldots
\nu_k}^*S_{\al_1\ldots \al_{i-1}\al_i\ldots \al_l}=$$ $$=S_{\nu_1\ldots
\nu_{i-1}}S_{\nu_i\ldots \nu_k}S_{\nu_i\ldots \nu_k}^*S_{\nu_1\ldots
\nu_{i-1}}^*S_{\nu_1\ldots \nu_{i-1}}S_{\al_i\ldots \al_l}=$$ $$=S_{\nu_1\ldots
\nu_{i-1}}S_{\nu_1\ldots \nu_{i-1}}^*S_{\nu_1\ldots \nu_{i-1}}S_{\nu_i\ldots
\nu_k}S_{\nu_i\ldots \nu_k}^*S_{\al_i\ldots \al_l}=0$$ because
$S_{\nu_i}^*S_{\al_i}=0$.
\end{enumerate}
\fim

\begin{teorema}
\label{t1} If $\nu\in P\backslash\La^*$ then $S_\nu=0$.
\end{teorema}
\demo Write $\nu=g_j\al$, and in this way,
$$S_\nu^*S_\nu=S_\alpha^*S_{g_j}^*S_{g_j}S_\alpha=S_\alpha
^*S_\alpha-\sum\limits_{k=1}^\infty\sum\limits_{\mu\in L_{g_j}^k
}S_\alpha^*S_\mu S_\mu^* S_\alpha.$$ We will analyse the summands of
$\sum\limits_{k=1}^\infty\sum\limits_{\mu\in L_{g_j}^k}S_\alpha^*S_\mu S_\mu^*
S_\alpha$ in the following way:

\begin{enumerate}
\item [] {\it Case 1: $|\mu|>|\alpha|$}\nl By Lemma
$\ref{lt3}$, $S_\mu S_\mu^*S_\al\neq 0$ only if $\mu=\alpha r$, for some $r\in
P$. We will show that there exists no such $r$. Suppose $\mu\in L_{g_j}^k$ is
such that $\mu=\al r$, with $|r|=l$. By definition of $L_{g_j}^k$,
$g_j\mu_1\ldots \mu_{k-1}\in \La^*$, but  $g_j\mu_1\ldots \mu_{k-1}=g_j\al
r_1\ldots r_{l-1}$, and so $\nu=g_j\al \in \La^*$. This is a contradiction,
because we are supposing $\nu\notin\La^*$. Therefore $\mu\neq \alpha r$, for
all $r\in P$, and so, by Lemma $\ref{lt3}$,\lb  $S_\al^*S_\mu
S_\mu^*S_\al=S_\al^*(S_\mu S_\mu^*S_\al)=0$ for all $\mu$ with $|\mu|>|\al|$.
\item [] {\it Case 2: $|\mu|\leq|\alpha|$}\nl By Lemma $\ref{lt3}$,
$S_\mu S_\mu^*S_\al\neq 0$, only if $\al=\mu r$, for some $r$ em $P$, and by
Lemma $\ref{lt1}$ if there exists such $\mu\in \cup L_{g_j}^k$, it is unique.
In this case we have by Lemma $\ref{lt3}$ that $S_\al^*S_\mu S_\mu^*
S_\al=S_\al^*(S_\mu S_\mu^* S_\al)=S_\al^*S_\al.$
\end{enumerate}
In this way, $S_\nu^*S_\nu=zS_\al^*S_\al$, where $z=0$ if there exists $\mu\in
\bigcup\limits_{k\in \N}L_{g_j}^k$ such that $\al=\mu r$ for some $r\in P$, and
$z=1$ otherwise.\nl Write $\nu=\nu_1\ldots \nu_k$ and take the smallest index
$i$ such that $\nu_{i+1}\ldots  \nu_{k}\in \La^*$. So,
$$S_\nu^*S_\nu=z_1S_{\nu_2\ldots \nu_k}^*S_{\nu_2\ldots \nu_k}=\ldots
=z_1\ldots  z_{i-1}S_{\nu_i\ldots \nu_k}^*S_{\nu_i\ldots \nu_k},$$ where $z_i$
are 0 or 1. We will show that $S_{\nu_i\ldots \nu_k}^*S_{\nu_i\ldots \nu_k}=0$.
Since $\nu_i\ldots \nu_k\notin\La^*$, by case 1 and case 2 above, we need to
show that there exist some $\mu\in \bigcup\limits_{k\in \N}L_{\nu_i}^k$ such
that $\nu_{i+1}\ldots \nu_k=\mu r$ for some $r\in P$.\nl Take the index $j$
such that $\nu_i\ldots \nu_j\in \La^*$ but $\nu_i\ldots
\nu_j\nu_{j+1}\notin\La^*$. Such index exists because $\nu_i\in\La^*$ and
$\nu_i\ldots\nu_k\notin\La^*$. Moreover, $\nu_{i+1}\ldots \nu_{j+1}\in \La^*$
because \lb $\nu_{i+1}\ldots \nu_k\in \La^*$, and so, $\nu_{i+1}\ldots \nu_{j+1}\in
L_{\nu_i}^{j+1-i}$. Thereby $S_{\nu_i\ldots \nu_k}^*S_{\nu_i\ldots \nu_k}=0$,
and so $S_{\nu}^*S_{\nu}=0$, in other words, $S_\nu=0$. \nl \fim

Observe that if $r=r_1\ldots r_k$ is in reduced form, with $r_i\in G^{-1}$ and
$r_{i+1}\in G$, then $S(r_ir_{i+1})=S(r_i)S(r_{i+1})=0$, from where $S(r)=0$.
Also, if $r=r_1\ldots r_k$ and $s=s_1\ldots s_l$ are elements of $\F_n$ in
reduced form and $r_k\neq s_1^{-1}$, then the reduced form of $rs$ is
$r_1\ldots r_ks_1\ldots s_l$, and so $S(rs)=S(r)S(s)$ by definition of $S$. \nl

\begin{definicao}
Given a group $\G$ and a Hilbert space $H$, a map $S:\G\rightarrow
\mathcal{L}(H)$ is a partial representation of the group $\G$ on $H$ if:
\begin{enumerate}
\item [$P_1$)] $S(e)=I$, where $e$ is the neutral element of $\G$ and $I$ is the identity operator on
$H$,
\item [$P_2$)] $S(t^{-1})=S(t)^*$, $\forall t\in \G$,
\item [$P_3$)] $S(t)S(r)S(r^{-1})=S(tr)S(r^{-1})$, $\forall t,r\in \G$,
\end{enumerate}
\end{definicao}

\begin{teorema}\label{t2}
If the map $S:\F_n\rightarrow \mathcal{L}(H)$ defined before satisfies
$M_1$,$M_2$ and $M_3$, then $S$ is a partial representation of the group $\F_n$
on $H$.
\end{teorema}
\demo Property $P_1$ is trivial. The proof of $P_2$  will be by induction on
$|t|$. If $|t|=1$, the equality between $S(t^{-1})$ and $S(t^*)$ is obviously
true. Suppose $S(t^{-1})=S(t^*)$ for all $t\in \F_n$ with $|t|=k$. Take
$t\in\F_n$ with $|t|=k+1$ and write $t=\til{t}x$, where $|\til{t}|=k$. Using
the induction hypothesis and the fact that the equality is true for $|x|=1$,
$$S(t^{-1})=S((\til{t}x)^{-1})=S(x^{-1}\til{t}^{-1})=S(x^{-1})S(\til{t}^{-1})=$$
$$S(x)^*S(\til{t})^*=(S(\til{t})S(x))^*=S(\til{t}x)^*=S(t)^*.$$ To verify
property $P_3$ we will prove the following :\nl{\it Claim: For all $r$ in
$\F_n$ and $t$ in $G\cup G^{-1}$, $E(r)=S(r)S(r)^*$ and $E(t)=S(t)S(t)^*$
commute.}\nl If $r=r_1\ldots  r_k$ where $r$ is in its reduced form, with
$r_i\in G^{-1}$ and $r_{i+1}\in G$ for some $i$, then $S(r)=0$ and so the claim
is trivial . Therefore let $r=\al \beta^{-1}$, where $r$ is in reduced form and
$\al, \beta\in P$. If $\beta\notin\La^*$, by Theorem $\ref{t1}$, $S_\beta=0$
from where we again see that the claim follows. Thus let us consider $\beta\in
\La^*$. \nl {\it Case 1: If $t\in G$, that is, $t=g_j$, for some $j$.}
\begin{enumerate}
\item [a)] $|\al|\neq 0$.\nl
Write $\al=\al_1\ldots \al_l$. If $\al_1\neq g_j$, then $S(g_j)^*S(\al)=0$ and
so \lb $E(t)E(r)=0=E(r)E(t)$. If $\al_1=g_j$ we have
$$S(\al)^*S(g_j)S(g_j)^*=S(\al_2\ldots \al_l)^*S(\al_1)^*S(g_j)S(g_j)^*=$$
$$=S(\al_2\ldots \al_l)^*S(\al_1)^*S(\al_1)S(\al_1)^*= S(\al_2\ldots
\al_l)^*S(\al_1)^*=$$ $$=(S(\al_1)S(\al_2\ldots \al_k))^*=S(\al)^*$$ and
similarly $S(g_j)S(g_j)^*S(\al)=S(\al)$. It follows that $E(t)$ and $E(r)$
commute.

\item [b)] $|\al|=0$.\nl
We have $r=\beta^{-1}$. Since $\beta\in\La^*$, using $M_2$,
$$E(r)E(t)=S(r)S(r)^*S(t)S(t)^*=S(\beta)^*S(\beta)S(g_j)S(g_j)^*=$$
$$=S(g_j)S(g_j)^*S(\beta)^*S(\beta)=S(t)S(t)^*S(r)S(r)^*=E(t)E(r).$$

\end{enumerate}
{\it Case 2: If $t\in G^{-1}$, namely, $t=g_j^{-1}$, with $g_j\in G$.}

Note that $$E(r)E(t)=E(r)S_tS_t^*=E(r)S_{g_j}^*S_{g_j}=
E(r)\left(I-\sum\limits_{k=1}^\infty\sum\limits_{\mu\in L_{g_j}^k}S_\mu
S_\mu^*\right)=$$ $$=E(r)-E(r)\left(\sum\limits_{k=1}^\infty\sum\limits_{\mu\in
L_{g_j}^k}S_\mu S_\mu^*\right)$$ and similarly,
$$E(t)E(r)=S_{g_j}^*S_{g_j}E(r)=E(r)-\left(\sum\limits_{k=1}^\infty\sum\limits_{\mu\in
L_{g_j}^k}S_\mu S_\mu^*\right)E(r).$$ To prove that $E(t)$ and $E(r)$ commute,
it is enough to show that $$E(r)S_\mu S_\mu^*=S_\mu
S_\mu^*E(r)\esp\esp\forall\mu\in L_{g_j}^k,\esp\esp \forall k\in \N.$$

\begin{enumerate}

\item [a)] $|\al|\neq 0$.
\begin{enumerate}
\item [i)] $|\al|\geq |\mu|$.\nl
By Lemma $\ref{lt3}$, if $\al=\mu s$ for some $s$ in $P$ then $S_\al^*S_\mu
S_\mu^*=S_\al^*.$ Therefore, $$E(r)S_\mu S_\mu^*=S_\al S_\beta^*S_\beta
S_\al^*S_\mu S_\mu^*=S_\al S_\beta^*S_\beta S_\al^*=E(r),$$ and similarly
$S_\mu S_\mu^*E(r)=E(r),$ and this proves that $E(r)S_\mu S_\mu^*=S_\mu
S_\mu^*E(r)$. Also by Lemma $\ref{lt3}$, if $\al\neq\mu s$ for all $s\in P$,
then $S_\al^*S_\mu S_\mu^*=0=S_\mu S_\mu^*S_\al$ and also in this case $E(r)$
and $S_\mu S_\mu^*$ commute.

\item [ii)] $|\al|<|\mu|$.\nl
By Lemma $\ref{lt3}$, if $\mu\neq \al s\esp\forall s\in P$, then $S_\al^*S_\mu
S_\mu^*=0=S_\mu S_\mu^*S_\al$, from where the equality follows. If $\mu=\al s$
for some $s\in P$, also by Lemma $\ref{lt3}$, $S_\al^*S_\mu S_\mu^*=S_sS_\mu^*$
and $S_\mu S_\mu^* S_\al=S_\mu S_s^*$, from where $$E(r)S_\mu S_\mu^*=S_\al
S_\beta^* S_\beta S_\al^*S_\mu S_\mu^*=S_\al S_\beta^* S_\beta S_s
S_\mu^*=S_\al S_\beta^* S_\beta S_s S_s^*S_\al^*,$$ and $$S_\mu
S_\mu^*E(r)=S_\mu S_\mu^* S_\al S_\beta^* S_\beta S_\al^*=S_\mu S_s^*S_\beta^*
S_\beta S_\al^*=S_\al S_sS_s^*S_\beta^* S_\beta S_\al^*.$$ Since $\beta\in
\La^*$, by $M_2$, $$S_s S_s^*S_\beta^* S_\beta=S_\beta^* S_\beta S_s S_s^*,$$
and this shows that $E(r)S_\mu S_\mu^*=S_\mu S_\mu^*E(r).$

\end{enumerate}

\item [b)] $|\al|=0$\nl
Since $\beta\in\La^*$, the equality between $E(r)S_\mu S_\mu^*$ and $S_\mu
S_\mu^*E(r)$ follows from $M_2$.

\end{enumerate}

This proves our claim. Let us now return to the proof of $P_3$, that is,
$$S(t)S(r)S(r^{-1})=S(tr)S(r^{-1}), \forall t,r\in \F_n.$$ To do this we use
induction on $|t|+|r|$. The equality is obvious if $|t|+|r|=1$. Suppose the
equality true for all $t,r\in \F_n$ such that $|t|+|r|<k$. Take $t,r\in \F_n$,
with $|t|+|r|=k$, write $t=\til{t}x, r=y\til{r}$, with $x,y\in G\cup G^{-1}$.
If $y\neq x^{-1}$, we have $S(tr)=S(t)S(r)$, from where
$S(tr)S(r^{-1})=S(t)S(r)S(r^{-1})$. Let us consider the case $x=y^{-1}$.
$$S(t)S(r)S(r^{-1})=S(\til{t}x)S(y\til{r})S((y\til{r})^{-1})=$$
$$=S(\til{t})S(x)S(y)S(\til{r})S(\til{r}^{-1})S(y^{-1})=S(\til{t})S(x)S(x^{-1})S(\til{r})S(\til{r}^{-1})S(x).$$
Using the claim and the fact that $S(x)$ is a partial isometry,
$$S(\til{t})S(x)S(x^{-1})S(\til{r})S(\til{r}^{-1})S(x)=S(\til{t})S(\til{r})S(\til{r}^{-1})S(x)S(x^{-1})S(x)=$$
$$=S(\til{t})S(\til{r})S(\til{r}^{-1})S(x)$$ and by the induction hypothesis,
$$S(\til{t})S(\til{r})S(\til{r}^{-1})S(x)=S(\til{t}\til{r})S(\til{r}^{-1})S(x).$$
On the other hand, $$S(tr)S(r^{-1})=S(\til{t}xy\til{r})S((y\til{r})^{-1})=$$
$$=S(\til{t}\til{r})S(\til{r}^{-1}y^{-1})=S(\til{t}\til{r})S(\til{r}^{-1})S(x).$$
This concludes the proof of $P_3$, and also of the theorem .\nl \fim

Instituto de Matem\'atica, Estat\'{\i}stica e Computa\c c\~ao Cient\'{\i}fica\\
UNICAMP-IMECC\\
Caixa Postal: 6065\\
13083-970 Campinas, SP, Brasil.\\
E-mail adress: royer@mtm.ufsc.br

\end{document}